# Eigenvalue Analysis via Kernel Density Estimation

Ahmed F. Yehia [a], Mohamed Saleh [a]

a Operations Research & Decision Support Department, Faculty of Computers and Information, Cairo University

**Abstract:**

In this paper, we propose an eigenvalue analysis -- of system dynamics models -- based on the Mutual Information measure, which in turn will be estimated via the Kernel Density Estimation method. We postulate that the proposed approach represents a novel and efficient multivariate eigenvalue sensitivity analysis.

## 1. Introduction and theoretical background:

The main contribution of this research is applying KDE on system dynamics models to measure the sensitivity of their eigenvalues to the input variables. This section is divided into the following subsections, in section 1.1. we briefly review eigenvalue analysis, while in section 1.2. we present the roadmap of the paper.

### 1.1. Eigenvalue Analysis:

Eigenvalue analysis studies the sensitivity of behavior modes (eigenvalues) – or their associated weights – to the gains of feedback loops. In general, the gains of feedback loops are functions of the gains of links; which in turn are functions of parameters, in the model. Hence one can link parameters (input variables) to certain eigenvalues (output variables). Eigenvalues are computed from the Jacobian matrix **J**, which can be calculated as illustrated in equation 1. After that, a sensitivity matrix $S_i$ is constructed for each eigenvalue, as illustrated in equation 2.

$$J = \begin{bmatrix} \partial \dot{X}_1/\partial X_1 & . & . & \partial \dot{X}_n/\partial X_1 \\ . & . & . & . \\ . & . & . & . \\ \partial \dot{X}_1/\partial X_m & . & . & \partial \dot{X}_n/\partial X_m \end{bmatrix} \quad (1)$$

$$S_i = \begin{bmatrix} \partial\lambda_i/\partial J(1,1) & . & . & \partial\lambda_i/\partial J(1,n) \\ . & . & . & . \\ . & . & . & . \\ \partial\lambda_i/\partial J(n,1) & . & . & \partial\lambda_i/\partial J(n,n) \end{bmatrix} \quad (2)$$

The majority of previous works in this area (Forrester 1983; Güneralp 2006; Kampmann and Oliva 2008; Gonçalves 2009; Saleh, Oliva et al. 2010), utilized the traditional univariate sensitivity (or elasticity) measure (explained above), which is based on the first order partial derivative of the eigenvalue with respect to an independent variable. This univariate measure represents only the marginal contribution of the independent variable – assuming that all other independent variables are constants. However, as Sterman put it: "*In nonlinear systems, the sensitivity of a system to variations in multiple parameters is not a simple combination of the response to the parameters varied alone*" (Sterman 2000).

Some scholars adopted multivariate linear regression analysis. Regression analysis is under the umbrella of multivariate analysis; i.e. can study simultaneous changes in more than one parameter (Esbensen, Guyot et al. 2002). Linear regression analysis aims to find a linear relationship between the dependent variable and independent variables. The core algorithm of linear regression is the least squares errors algorithm, which is used in data fitting (Rencher and Christensen 2012). Yet, the main limitation here is the assumption of a linear relationship between the dependent variable and independent variables.

Nonlinear regression analysis can be used when linear regression fails. However, the analyst must try different functions; e.g. exponential, logarithmic, etc. In addition initial values for the coefficients are needed; and in general, the solution changes according to initial settings (Draper, Smith et al. 1966; Hair, Black et al. 2006). That is, there is no guarantee to reach the global minimum (least squares errors); as there might be several local minima. Moreover, the solution of complex nonlinear regression equations might not converge. Finally, nonlinear multivariate regression is not well suited to handle non-monotonic functions (Rencher and Christensen 2012).

To summarize, traditional eigenvalue analysis methods have the following shortcomings or a subset of them, first, they measure only the marginal contribution of the independent variable, second, they capture only linear relationships and finally, they only perform univariate sensitivity analysis.

## 1.2. Roadmap of the Paper

The rest of this paper is organized as follows, in section 2, we give a brief theoretical background. In section 3, our proposed method is introduced. In section 4, we test our proposed method on a hypothetical model. Finally, we conclude and propose future work in section 5.

## 2. Theoretical Background:

Previously (Yehia, Saleh et al. 2014; Yehia, Saleh et al. 2015), we have applied Maximal Information Coefficient (MIC) and Kernel Canonical Correlation Analysis (KCCA) to eigenvalue analysis -- each of which had some limitations. MIC can only represent many-to-one relationships and has a very high computational time because of its recursive calls; on the other hand, KCCA has the extra complexity of choosing the appropriate kernel. Therefore, we continued our investigation for a better method to measure eigenvalue sensitivity. In this paper, our proposed method is based on Kernel Density Estimation (KDE) and Mutual Information (MI).

To explain the concepts related to the proposed method, this section is divided into the following subsections: section 2.1 briefly explains the Mutual Information (MI) concept; while section 2.2 explains the Kernel Density Estimation (KDE) concept.

### 2.1. Mutual Information (MI):

As stated before (Yehia, Saleh et al. 2014; Yehia, Saleh et al. 2015), we used a Maximal Information Coefficient (MIC) that is based on Mutual Information concept (Reshef, Reshef et al. 2011). Although MIC has advanced capabilities and can detect all kinds of relationships. It has some limitations regarding applying it to eigenvalue analysis, such as it only represents many-to-one relationships and has a very high computational time because of its recursive calls.

In this paper, we will overcome the limitations of the MIC, as we will focus on using a Multivariate Mutual Information to measure the total correlation between the SD model inputs and its eigenvalues. Although MI doesn't assume any prior functional relationship, it can detect any type of complex relationships (Moon, Rajagopalan et al. 1995; Steuer, Kurths et al. 2002; Dionisio, Menezes et al. 2004); and that is our main motivation to use it. In the rest of this subsection we will give a brief review about MI, but as an earlier step, we will review the Information Entropy (IE) which is the building block of the MI, as shown in the following equation. Note that the information entropy can be associated with a single univariate distribution (of a single variable), or associated with a joint distribution (of several variables).

$$MI(x,y) = IE(x) + IE(y) - IE(x,y) \qquad (3)$$

Information Entropy (IE) is a measure of the uncertainty of a probability distribution associated with a specific random variable (or set of variables). It's score ranges from 0 to 1, and it can be interpreted as the expected surprisal associated with the probability distribution. Note that the surprisal of an event is inversely proportional to its probability of occurrence. Information Entropy can be calculated using the following equation (Anderson 2008).

$$IE(x) = -\sum_{i=0}^{n} P(x_i) * log_2 p(x_i) \qquad (4)$$

On the other hand, Mutual Information measures how much information one variable can tell about another variable on average. Also one might see it is a dimensionless quantity that measures the reduction in uncertainty (info. entropy) about one random variable given knowledge of another. The following equations are the formal definitions of mutual information in the continuous and the discrete cases respectively (Gallager 1968).

$$MI(X;Y) = \int \int p(X,Y) \, log_2 \frac{p(X,Y)}{p(X)\,p(Y)} \, dx \, dy \qquad (5)$$

$$MI(x,y) = \sum_x \sum_y P(x,y) log \frac{P(x,y)}{P(x)P(y)} \qquad (6)$$

In the discrete case Mutual Information score should be normalized using the following equation in order to obtain a score between the 0 and 1, which is the same equation used in maximal information coefficient (MIC). The normalization is calculated simply based on the number of bins (n).

As one might notice, the Mutual Information equation is bivariate which is not consistent with our objective in this paper. So, in order to perform a multivariate sensitivity analysis, we utilized the generalized MI equation which is illustrated below (Hamming 1986).

$$MI(x_1, x_2, \ldots, x_m, y_1, y_2, \ldots, y_n) =$$

$$\sum_{x_1} \sum_{x_2} \ldots \sum_{x_m} \sum_{y_1} \sum_{y_2} \ldots \sum_{y_n} P(x_1, x_2, \ldots, x_m, y_1, y_2, \ldots, y_n) log \frac{P(x_1,x_2,\ldots,x_m,y_1,y_2,\ldots,y_n)}{P(x_1)*P(x_2)*\ldots*P(x_m)*P(y_1)*P(y_2)*\ldots*P(y_n)} \qquad (7)$$

The above equation can be expressed in terms of marginal entropies and joint entropy as follows (Hamming 1986):

$$MI(x_1, x_2,.., x_m, y_1, y_2, ..., y_n) = \sum_{i=1}^{m} IE(x_i) + \sum_{j=1}^{n} IE(y_j) - IE(x_1, x_2,.., x_m, y_1, y_2, ..., y_n) \quad (8)$$

In the above equation, in order to able to calculate the Information Entropy and Mutual Information we need to compute joint and marginal distributions. This can be done via Kernel Density Estimation (KDE) method explained in the next section.

## 2.2. Kernel Density Estimation (KDE):

Density Estimation methods belong to the nonparametric inferential statistics class. Inferential statistics is the class of statistics which can be used to infer information about the population using the available sample data points of the population. On the other hand, nonparametric statistics is the class of statistics that doesn't require a prior knowledge about the population distribution or parameters (Gibbons 1993; Asadoorian and Kantarelis 2005).

There are different density estimators, the most popular one and the easiest one to interpret is the histogram density estimator; however, it has several limitations. Moreover other density estimators – like kernel density estimation – are known to provide a lower mean squared error between the estimated and empirical cumulative distribution (Silverman 1986; Terrell 1990; Moon, Rajagopalan et al. 1995).

Kernel Density Estimation infer or estimate the probability density distribution of the population using the available finite population sample via different kernel functions. The estimation is the sum of the kernels associated with the data points. The density of the estimated distribution at a certain point is a translation of how many kernels intersect at this point. The population density distribution can be estimated using the following function (Silverman 1986):

$$\hat{f}(x, h) = \frac{1}{nh} \sum_{i=0}^{n} K\left(\frac{x - x_i}{h}\right) \quad (9)$$

The cumulative distribution of the above probability distribution function can be calculated as follows:

$$\hat{F}(x) = \int_{-\infty}^{x} \hat{f}(x, h) \, dx \quad (10)$$

Where $K(.)$ is the kernel function and $h$ is the density distribution smoothing parameter or the kernel bandwidth. The estimation process is an optimization problem of choosing the bandwidth that minimizes the mean squared error between the empirical cumulative distribution function and the estimated one.

The optimal kernel bandwidth $\boldsymbol{h}$ can be estimated via an optimization problem that seeks to minimize the mean squared error between the empirical cumulative distribution and the above estimated cumulative distribution.

There are several kernel functions that can be used. In this paper the Gaussian kernel (the most commonly used kernel function) is been adopted, due to its flexibility in estimating a wide variety of distributions. It can be calculated using the following equation (Silverman 1986).

$$K(y) = \frac{1}{\sqrt{2\pi}} e^{\frac{-y^2}{2}} \qquad (11)$$

Where $y = \dfrac{x - x_i}{h}$

Conceptually, a Gaussian/normal distribution is associated with each data point. Then the small distributions are summed up to estimate the total density function.

In our work, in order to conduct a multivariate sensitivity analysis, we need to estimate joint probability distribution associated with more than one variable. After that we could use the estimated joint distribution to calculate the marginal distribution for each variable. In general, multivariate kernel density estimation can be utilized for more than two variables. There is an analogy between the univariate and the multivariate case, the only difference between them is the calculation of the kernel bandwidth. In the univariate case, the bandwidth is a scalar variable $\boldsymbol{h}$; while in the multivariate case, the bandwidth is a matrix $\boldsymbol{H}$. It can be calculated as indicated in the following equation (Scott 2015).
.

$$\hat{f}(x)_H = \frac{1}{n} \sum_{i=0}^{n} K_H(X - X_i) \qquad (12)$$

Recall that, in the univariate case one can formulate an optimization problem to estimate the kernel bandwidth scalar variable. In the multivariate case a similar optimization problem can be formulated, but in this case the elements of the bandwidth matrix will be the decision variables of optimization problem.

While the objective function remains the same (i.e. minimize the mean squared error between the empirical CDF and the estimated CDF).

## 3. Proposed Method:

Our proposed method is based on KDE and MI; this section demonstrates the algorithm of the method which consists of six main steps listed below:

**First Step:** Normalize the variables between 0 and 1

**Second Step:** Estimate joint CDF in $R^n$ via KDE.

- *KDE is calculated by **Statsmodels** (python library).*

**Third Step:** Discretize the estimated joint CDF:

   **3.1.** Divide $R^n$ space into small N-Dimensional hypercubes
   **3.2.** Compute the probability in each cube from the estimated joint CDF.

**Fourth Step:** Compute the marginal probability distributions.

**Fifth Step:** Compute marginal and joint entropies.

**Sixth Step:** Compute Mutual Information.

**Seventh Step:** Compute the sensitivities of the output variables.

The rest of this section is an explanation of the above steps. The python codes for the above steps are documented in appendix A.

In the first step, all the model input variables and eigenvalues are normalized between zero and one.

In the second step, Gaussian kernel density estimation is utilized to estimate the joint cumulative distribution function for all the variables combined.

The third step is divided into two main sub-steps, in the first one, we divide the $R^n$ space into n-dimensional hypercubes using equally sized bins. In the second one, we compute the probability in each hypercube using the estimated joint cumulative distribution function.

In the fourth step, we compute the marginal distribution associated with each variable (input variables and eigenvalues). This is done by summing the probabilities in all dimensions except for the intended variable.

In the fifth step, the joint and marginal entropies are calculated using the estimated joint CDF and the marginal distributions calculated in step four.

In the sixth step, we compute Mutual Information using the generalized formula specified in equation 8.

In the last step, we compute the sensitivity associated with each input variables. Sensitivity of the model eigenvalues (output) to a specified input variable will be the reduction in the Mutual Information score caused by eliminating that. The higher the amount reduced from the MI score; the higher influence the eliminated variable and the opposite is true.

In the next section, we will test the algorithm on a hypothetical model. Appendix B documents the entire python program for conducting this experiment.

## 4. Experiments:

This section shows the experiments conducted in order to apply the KDE and MI to perform eigenvalue analysis on a simple hypothetical system dynamics model. The stock and flow diagram of the model is shown in the following figure. The model consists of two stocks: $S_1$ and $S_2$. While $R_1$ and $R_2$ are the inflows of $S_1$ and $S_2$ respectively. Moreover, there are three auxiliary variables: $g_{11}$, $g_{12}$, $g_{21}$. The equations of the model are as follows:

- $X_1$, $X_2$ & $X_3$ are uniform random variables [0, 1]. These are the input variables (parameters)
- $g_{11} = X_1 * X_2 * X_3$
- $g_{12} = X_1 * X_2$
- $g_{21} = X_1 * X_3$

- $R_1 = g_{11}*S_1 + g_{12}*S_2 + 1$
- $R_2 = S_1*g_{21}$

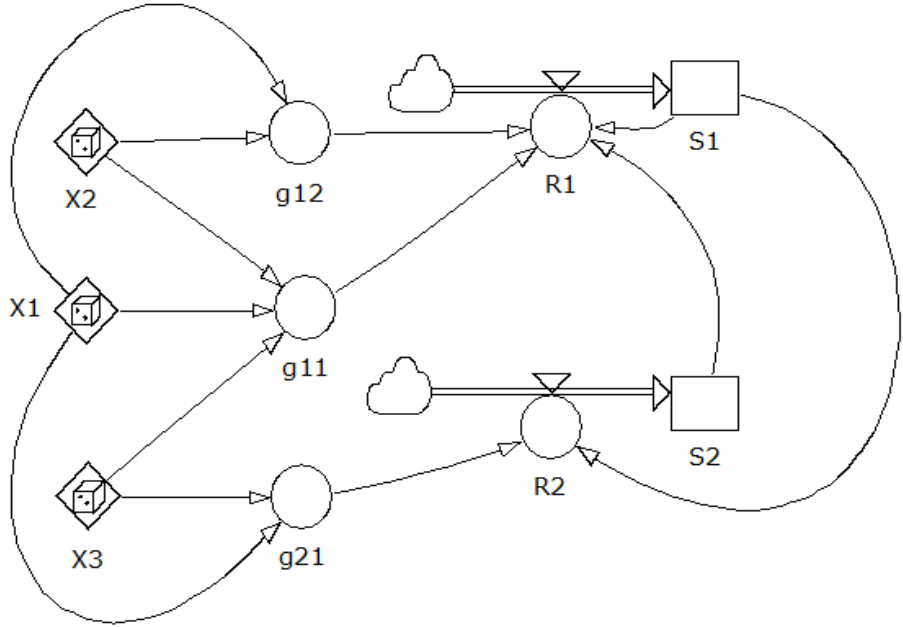

*Figure 1: The Stock & Flow Diagram of the Model used in the Eigenvalue Experiments.*

Referring to section 1.2, model eigenvalues can be computed from the model Jacobian matrix (which can be considered a condensed mathematical representation of the model structure). Our model has two eigenvalues as it has two stocks only; moreover, these eigenvalues are always a real number (i.e. not complex numbers) – regardless the values of the parameters. The results of the experiments are shown in the following two tables. Table 1, shows the results associated with a single output variable, which is the dominant eigenvalue; in our case y2 is strictly dominating y1, so we will refer to y2 as the dominant eigenvalue or Ymax. While, table 2 shows the results associated with two output variables, which are the two eigenvalues of the model. All the experiments are conducted under the following setup:

1. Sample size = 12000
2. #cubes per dimension = 10

*Table 1: The association level between all sets of input variables and the dominant eigenvalue (Ymax)*

| Ymax | X1 | X2 | X3 | X1 & X2 | X1 & X3 | X2 & X3 | X1, X2 & X3 |
|---|---|---|---|---|---|---|---|
| **MI** | 0.15 | 0.1 | 0.09 | 0.28 | 0.27 | 0.21 | 0.5 |

*Table 2: The association level between all sets of input variables and the model eigenvalues*

| Y1 & Y2 | X1 | X2 | X3 | X1 & X2 | X1 & X3 | X2 & X3 | X1, X2 & X3 |
|---|---|---|---|---|---|---|---|
| **MI** | 0.65 | 0.59 | 0.58 | 0.78 | 0.77 | 0.7 | 0.99 |

By inspecting the previous two tables, one can notice that; the mutual information between all the model parameters and the two eigenvalues equals to one (as it represents a full information set). Also, the Mutual Information between all the model input variables and one of the eigenvalues equals 0.5 -- which exactly equals to the mutual information between the two eigenvalues. I.e.

$$MI\ (y1,\ y2) = 0.5$$

In Addition, both tables provide indications about the significance of the input variables. For example, X1 is more significant than X2 and X3.

The bottom line of this method is how good the estimated joint distribution represents the distribution of actual sample points and hence the population. This optimal representation can be obtained by minimizing the least square error between the sample distribution and the estimated distribution -- via searching for the optimal the kernel bandwidth matrix. Moreover, it is possible to use Silverman (Silverman 1986) rule of thumb, that provides an easy heuristic (which skips the optimization process) to set the bandwidth size using an approximation formula.

# 5. Conclusion:

This paper presents a multivariate nonlinear eigenvalue analysis method for dynamic model. The method is based on combining Kernel Density Estimation and Mutual Information. The proposed method overcomes the research gaps of the traditional methods. Specifically, the proposed method is a multivariate (many-to-many) sensitivity analysis methods, which can detect highly nonlinear relationships between model features and its eigenvalues. Also, it has a reasonable computational time especially using the Silverman's approximation formula. Moreover, there is no need to simulate the underline model. The proposed solution only interacts with a condensed matrix representation of the model; i.e. the Jacobian matrix. In each run, the process, which takes time, is the computation of eigenvalues from the Jacobian matrix; and there are algorithms that compute eigenvalues very fast. Via this method, decision-makers can rank policy parameters according to their impacts on the dominant eigenvalues. The downside of the proposed method is the limited number of variables, typically six variables are the practical upper limit for the kernel density estimation method.

In future, we will continue the experimental work, and test the method on more complex highly nonlinear models. In addition, we plan to develop a wrapper-based parameter selection method to facilitate the automatic ranking and selection of parameters. Finally, in many nonlinear dynamic models, eigenvalues depend on the current state of the model. For these cases, we plan to devise an innovative framework that links the time trajectory of the dominant eigenvalues with parameters, in the model.

# Appendix A: Step by Step Python Code for the Proposed Method

**Step Two:**

### Calculating the joint distribution using statsmodels library.

dens_fn = sm.nonparametric.KDEMultivariate(data=[x1, x2, x3, norm_y1, norm_y2], var_type='ccccc', bw='cv_ls')

**Step Three:**

### The two substeps of the third step.

def sum_prob_origion_cur(vert_ind_vect, probNCube, cur_vert_index):

   total = 0

   for i in range(0, len(vert_ind_vect)):

      if True not in np.greater_equal(cur_vert_index, vert_ind_vect[i]):

      break

     if False in np.greater_equal(cur_vert_index, vert_ind_vect[i]):

      continue

     total += probNCube[i]

   return total

def discreatize (dens_fn, n_dim, nCube_perDim):

   probNCube = np.zeros(int(nCube_perDim ** n_dim))

   cdfMatrix = np.zeros(int((nCube_perDim) ** n_dim)+1)

   NCube_len = 1.1/nCube_perDim

   vert_ind_vect = np.mgrid[[slice(0, nCube_perDim, 1) for n in range(0, n_dim)]].reshape(n_dim,-1).T

   for i in range(0, len(probNCube)):

     cdfMatrix[i+1] = dens_fn.cdf((vert_ind_vect[i]+1) * NCube_len) ;

     probNCube[i] = cdfMatrix[i+1] - sum_prob_origion_cur(vert_ind_vect, probNCube, vert_ind_vect[i]);

   return probNCube.reshape([nCube_perDim for n in range(0, n_dim)]).T

probNCube = discreatize(dens_fn, n_dim, nCube_perDim)

**Step Four:**

### Calculate the marginal distributions using the discretized estimated joint CDF.

m0 = np.sum(np.sum(np.sum(np.sum(probNCube,4) ,3) ,2), 1)

```
m1 = np.sum(np.sum(np.sum(np.sum(probNCube,4) ,3) ,2), 0)

m2 = np.sum(np.sum(np.sum(np.sum(probNCube,4) ,3) ,1), 0)

m3 = np.sum(np.sum(np.sum(np.sum(probNCube,4) ,2) ,1), 0)

m4 = np.sum(np.sum(np.sum(np.sum(probNCube,3) ,2) ,1), 0)
```

**Step Five:**

```
### Calculate marginal and joint entropies.
def entropy(X):
    X[X<=0.00001] = 1
    entropy = sum(p*np.log2(1/p) for p in X)
    X[X==1] = 0
    return entropy

def joint_entropy(probNCube, n_dim):
    joint_vector = probNCube.reshape(len(probNCube)**n_dim, 1)
    return entropy(joint_vector)[0]
```

**Step Six:**

```
### calculating Mutual Information
MI = entropy_x1 + entropy_x2 + entropy_x3 + entropy_y1 + entropy_y2 - entropy_joint
```

**Step Seven:**

```
### calculating another Mutual Information score without the first input variable (x1)
MI_without_x1 = entropy_x2 + entropy_x3 + entropy_y1 + entropy_y2 - entropy_joint

### calculate the sensitivity of the model eigenvalues to the first input variable.
Sensitivity(x1) = MI - MI_without_x1
### repeat for all input variable
```

# Appendix B: Python Code for the Proposed Method

```
import numpy as np
import statsmodels.api as sm
```

```python
def entropy(X):
    X[X<=0.00001] = 1
    entropy = sum(p*np.log2(1/p) for p in X)
    X[X==1] = 0
    return entropy

def joint_entropy(probNCube, n_dim):
    joint_vector = probNCube.reshape(len(probNCube)**n_dim, 1)
    return entropy(joint_vector)[0]

def sum_prob_origion_cur(vert_ind_vect, probNCube, cur_vert_index):
    total = 0
    for i in range(0, len(vert_ind_vect)):
            if True not in np.greater_equal(cur_vert_index, vert_ind_vect[i]):
            break
        if False in np.greater_equal(cur_vert_index, vert_ind_vect[i]):
            continue
        total += probNCube[i]
    return total

def discreatize (dens_fn, n_dim, nCube_perDim):
    probNCube = np.zeros(int(nCube_perDim ** n_dim))
    cdfMatrix = np.zeros(int((nCube_perDim) ** n_dim)+1)
    NCube_len = 1.1/nCube_perDim
    vert_ind_vect = np.mgrid[[slice(0, nCube_perDim, 1) for n in range(0, n_dim)]].reshape(n_dim,-1).T
    for i in range(0, len(probNCube)):
        cdfMatrix[i+1] = dens_fn.cdf((vert_ind_vect[i]+1) * NCube_len) ;
        probNCube[i] = cdfMatrix[i+1] - sum_prob_origion_cur(vert_ind_vect, probNCube, vert_ind_vect[i]);
    return probNCube.reshape([nCube_perDim for n in range(0, n_dim)]).T
```

```
### main
n_dim = 5
nCube_perDim = 10

dens_fn = sm.nonparametric.KDEMultivariate(data=[x1, x2, x3, norm_y1, norm_y2], var_type='ccccc', bw='cv_ls')

probNCube = discreatize(dens_fn, n_dim, nCube_perDim)

m0 = np.sum(np.sum(np.sum(np.sum(probNCube,4) ,3) ,2), 1)
m1 = np.sum(np.sum(np.sum(np.sum(probNCube,4) ,3) ,2), 0)
m2 = np.sum(np.sum(np.sum(np.sum(probNCube,4) ,3) ,1), 0)
m3 = np.sum(np.sum(np.sum(np.sum(probNCube,4) ,2) ,1), 0)
m4 = np.sum(np.sum(np.sum(np.sum(probNCube,3) ,2) ,1), 0)

entropy_x1 = entropy(m0)
entropy_x2 = entropy(m1)
entropy_x3 = entropy(m2)
entropy_y1 = entropy(m3)
entropy_y2 = entropy(m4)
entropy_joint = joint_entropy(probNCube, n_dim)

MI = entropy_x1 + entropy_x2 + entropy_x3 + entropy_y1 + entropy_y2 - entropy_joint
```

### example: calculating the sensitivity of the eigenvalues to x1, remove x1 from calculations.

dens_fn = sm.nonparametric.KDEMultivariate(data=[x2, x3, norm_y1, norm_y2], var_type='cccc', bw='cv_ls')

probNCube = discreatize(dens_fn, n_dim, nCube_perDim)

m0 = np.sum(np.sum(np.sum(np.sum(probNCube,3) ,2), 1)
m1 = np.sum(np.sum(np.sum(np.sum(probNCube,3) ,2), 0)
m2 = np.sum(np.sum(np.sum(np.sum(probNCube,3) ,1), 0)
m3 = np.sum(np.sum(np.sum(np.sum(probNCube,2) ,1), 0)

entropy_x2 = entropy(m0)
entropy_x3 = entropy(m1)
entropy_y1 = entropy(m2)
entropy_y2 = entropy(m3)
entropy_joint = joint_entropy(probNCube, n_dim)

MI_without_x1 = entropy_x2 + entropy_x3 + entropy_y1 + entropy_y2 - entropy_joint

sensitivity_x1 = MI − MI_without_x1